\documentclass[11pt]{article}

\usepackage{a4wide}
\usepackage{amssymb}
\usepackage{amsmath}
\usepackage{theorem}
\usepackage{epsfig}
\usepackage{verbatim}
\usepackage{graphicx}
\usepackage{color}
\usepackage{indentfirst}
\newtheorem{thm}{Theorem}[section]

\newtheorem{rem}[thm]{Remark}

\newtheorem{lemma}[thm]{Lemma}

\newcommand{\R}{\Bbb{R}}

\newcommand{\D}{\displaystyle}

\newcommand{\sign}{{\rm sign}\thinspace}

\newcommand{\lae}{\Lambda^{\alpha}}

\newcommand{\al}{\alpha}
\newcommand{\ep}{\varepsilon}

\newcommand{\om}{\Omega}
\newcommand{\ff}{\mathbb}

\newcommand{\pq}{_q^p}
\newcommand{\be}{\beta}
\newcommand{\pa}{\partial}

\def\eop{{\ \vrule height 3pt width 3pt depth 0pt}}

\begin{document}

\author{Angel Castro, Diego C\'ordoba and Francisco Gancedo}
\title{Singularity Formation in a Surface Wave Model}
\date{April 21, 2010}
\maketitle
\begin{abstract}

In this paper we study the Burgers equation with a nonlocal term of the form $Hu$ where $H$ is the Hilbert transform. This system has been considered as a quadratic approximation for the dynamics of a free boundary of a vortex patch (see \cite{marsden} and \cite{HB}). We prove blow up in finite time for a large class of initial data with finite energy. Considering a more general nonlocal term, of the form $\lae Hu$ for $0<\al< 1$, finite time singularity formation is also shown.

\end{abstract}
\section{Introduction.}
We shall study the formation of singularities for the equation
\begin{eqnarray}
u_t+uu_x&=&\Lambda^\alpha Hu, \label{BLH}\\
u(x,0)&=&u_0(x)\nonumber,
\end{eqnarray}
with $0\leq\alpha<1$, where $H$ is the Hilbert transform \cite{St1} defined by
$$Hf(x)=\frac{1}{\pi}P.V.\int_{\ff{R}}\frac{f(y)}{x-y}dy$$
and  $\Lambda^\al\equiv (-\Delta)^{\al/2}$ is given by the following expression
$$\Lambda^\al f(x)=k_\alpha \int_{\ff{R}}\frac{f(x)-f(y)}{|x-y|^{1+\al}}dy,\qquad  k_\al=\frac{\Gamma(1+\alpha)\cos((1-\alpha)\pi/2)}{\pi}.$$
The case $\alpha=0$
\begin{equation}
u_t+uu_x= Hu \label{BH}
\end{equation}
was introduced by J. Marsden and A. Weinstein \cite{marsden} as a second order approximation for the dynamics of a free boundary of a vortex patch (see \cite{Ch} and \cite{BC}). Recently J. Biello and J.K. Hunter \cite{HB} proposed it as a model for waves with constant nonzero linearized frequency. They gave a dimensional argument to show that it models nonlinear Hamiltonian waves with constant frequency. In addition, an asymptotic equation from (\ref{BH}) is derived, describing surface waves on a planar discontinuity in vorticity for a two-dimensional inviscid incompressible fluid. They also carried out  numerical analysis showing evidence of singularity formation in finite time. Let us point out that the Hamiltonian structure of the equation (\ref{BLH}) (in particular for $\al=0$) comes from the representation
\begin{eqnarray}\label{hamiltonian}
u_t+\pa_x\left[\frac{\delta \mathcal{H}}{\delta u}\right]=0,\quad&\mbox{where}&\quad \mathcal{H}(u)=\D\int_{\R}\Big(\frac{1}{2}u\Lambda^{\alpha-1}u+\frac{1}{6}u^3\Big)dx.
\end{eqnarray}

In section \ref{seccion1} we show that the linear term in the equation (\ref{BH}) is too weak to prevent the singularity formation of the Burgers equation. In fact, we show that, if the $L^\infty$ norm of the initial data is large enough compare with the $L^2$ norm, the maximum of the solution has a singular behavior  during the time of existence. One of the ingredients in the proof is to use the following pointwise inequality
  \begin{equation}\label{pwi}
  u(x)^4\leq 16||u||^2_{L^2(\ff{R})}\int_{\ff{R}}\frac{(u(x)-u(y))^2}{(x-y)^2}dy,
  \end{equation}
(see lemma \ref{cuarta} below) which can be understood as the local version of the well-known bound
$$
||u||^4_{L^4}\leq C||u||^2_{L^2}||\Lambda^{1/2} u||^2_{L^2}=\frac{C}{2\pi}||u||^2_{L^2}\int_{\ff{R}}\int_{\ff{R}}\frac{(u(x)-u(y))^2}{(x-y)^2}dydx.
$$
In the appendix we provide a generalized pointwise inequality ($n-$dimensional) in terms of fractional derivatives.

 In section 3 we consider the more general family of equations, with a higher order term in derivatives, given by (\ref{BLH}). By a different method, we  prove that the blow up phenomena still arises. Let us note that, since $\Lambda H u=-u_x$, the case $\al=1$ trivializes. Using the same approach as in section 2, it is possible to obtain blow up for $0<\al<1/3$. Inspired by the method used in \cite{ADM}, we check the evolution of the following quantity
$$J_q^p u(x)=\int_{\ff{R}}w\pq(x-y)u(y)dy,\quad \mbox{where}\quad w\pq(x)= \left\{\begin{array}{ccc}|x|^{-q}\sign(x) &\text{if} \, |x|<1\\ |x|^{-p} \sign(x)&\text{if}\, |x|>1\end{array}\right.,
$$
with $0<q<1$ and $p>2$ to find a singular behavior. Let us note that a similar approach was used by H. Dong, D. Lu and D. Li (see \cite{DDL}) to show blow up for the Burgers equation with fractional dissipation in the supercritical case ($0<\al<1$):
\begin{equation}\label{dibur}
u_t+ uu_x= -\lae u.
\end{equation}
A different method to show singularities can be found in \cite{KNS}.

It is well known that the $L^p$ norms of the solutions of equation (\ref{dibur}) are bounded for all $1\leq p\leq \infty$. However, to the best to the authors knowledge,  two quantities are conserved by equation (\ref{BLH}). The orthogonality property of the Hilbert transform provides the conservation of the $L^2$ norm, i.e.
\begin{equation*}
||u(\cdot,t)||_{L^2(\ff{R})}=||u_0||_{L^2(\ff{R})}.
\end{equation*}
Since the equation is given by (\ref{hamiltonian}), we have that
\begin{equation*}
\int_{\ff{R}}\Big(\frac{1}{3}u^3(x,t)+\Big(\Lambda^\frac{\al-1}{2}u(x,t)\Big)^2\Big)dx=\int_{\ff{R}}\Big(\frac{1}{3}u_0^3(x)+\Big(\Lambda^\frac{\al-1}{2}u_0(x)\Big)^2\Big)dx.
\end{equation*}

\section{Blow up for the Burgers-Hilbert equation.}\label{seccion1}

The purpose of this section is to show finite time singularity formation in solutions of the equation (\ref{BH}). The result we shall prove is the following:

\begin{thm}\label{blowup}

Let $u_0 \in L^2(\ff{R})\cap C^{1+\delta}(\ff{R})$, with $0<\delta<1$, satisfying the following condition:

There exists a point $\be_0\in \ff{R}$ with
\begin{equation}\label{uno}Hu_0(\be_0)>0,\end{equation}
such that
\begin{equation}\label{dos}u_0(\be_0)\geq \left(32\pi ||u_0||_{L^2(\ff{R})}^2\right)^{1/3}.\end{equation}
Then there is a finite time $T$ such that

$$\lim_{t\rightarrow T}||u(\cdot,t)||_{C^{1+\delta}(\ff{R})}=\infty,$$
where $u(x,t)$ is the solution to the equation (\ref{BH}).
\end{thm}

Proof: Let us assume that there exist a solution of the equation (\ref{BH}) $$u(x,t)\in C([0,T),C^{1+\delta}(\ff{R})),$$ for all time $T<\infty$ and with $u_0$ satisfying the hypotheses.

Now, we shall define the trajectories $x(\be,t)$  by the equation
\begin{eqnarray*}
\frac{d x(\be,t)}{dt}&=&u(x(\be,t),t),\\
x(\be,0)&=&\be.
\end{eqnarray*}
Considering the evolution of the solution along trajectories, it is easy to get the identity
$$\frac{du(x(\be,t),t)}{dt}=u_t(x(\be,t),t)+\frac{d x(\be,t)}{dt}u_x(x(\be,t),t)=Hu(x(\be,t),t),$$
and taking a derivative in time we obtain
\begin{align*}
\frac{d^2u(x(\be,t),t)}{dt^2}&=Hu_t(x(\be,t),t)+ u(x(\be,t),t) Hu_x(x(\be,t),t)\\
&=-H(uu_x)(x(\be,t),t)-u(x(\be,t),t)+u(x(\be,t),t)Hu_x(x(\be,t),t).
\end{align*}
Since
$$H(uu_x)(x)=\frac{1}{2}H((u^2)_x)=\frac{1}{2}\Lambda(u^2)(x),$$
we can write
$$\frac{1}{2}\Lambda(u^2)(x)=\frac{1}{2\pi}P.V\int_{\ff{R}}\frac{u(x)^2-u(y)^2}{(x-y)^2}dy=u(x)\Lambda u(x)-\frac{1}{2\pi}\int_{\ff{R}}\frac{(u(x)-u(y))^2}{(x-y)^2}dy,$$
 and therefore it follows that
\begin{equation}\frac{d^2u(x(\be,t),t)}{dt^2}=\frac{1}{2\pi}\int_{\ff{R}}\frac{(u(x(\be,t),t))-u(y,t)))^2}{(x(\be,t)-y)^2}dy-u(x(\be,t),t).\label{segunda}\end{equation}

In order to continue with the proof we will prove the lemma below (for similar approach see \cite{CC}):

\begin{lemma}\label{cuarta}
Let $u\in L^2(\ff{R})\cap C^{1+\delta}(\ff{R})$, for $0<\delta<1$. Then
$$\frac{1}{2\pi}\int_{\ff{R}}\frac{(u(x)-u(y))^2}{(x-y)^2}dy\geq C u(x)^4,$$
where
$$C=\frac{1}{32\pi E}$$
and
$$E=||u||_{L^2(\ff{R})}^2.$$

\end{lemma}

Proof of lemma \ref{cuarta}: Let us assume that $u(x)> 0$ (a similar proof holds for $u(x)<0$). Let $\Omega$ be the set
$$\Omega=\{y\in \ff{R}\quad:\quad |x-y|<\Delta\},$$
where $\Delta$ will be given below. And let $\Omega^1$ and $\Omega^2$ be the subsets
\begin{eqnarray*}
\Omega^1&=&\{y\in \Omega\, : u(x)-u(y)\geq\frac{u(x)}{2}\},\\
\Omega^2&=&\{y\in \Omega\, : u(x)-u(y)<\frac{u(x)}{2}\}=\{y\in \Omega \, : u(y)>\frac{u(x)}{2}\}.
\end{eqnarray*}
Then
$$\frac{1}{2\pi}\int_{\ff{R}}\frac{(u(x)-u(y))^2}{(x-y)^2}dy\geq \frac{u(x)^2}{8\pi\Delta^2} |\Omega^1|.$$
On the other hand
$$E=\int_{\ff{R}}u(y)^2dy\geq \int_{\Omega^2}u(y)^2dy\geq \frac{u(x)^2}{4}|\Omega^2|,$$
and therefore
$$|\Omega^2|\leq\frac{4E}{u(x)^2}.$$
Since $|\Omega^1|=|\Omega|-|\Omega^2|$ and $|\Omega|=2\Delta$, we have that
$$\frac{1}{2\pi}\int_{\ff{R}}\frac{(u(x)-u(y))^2}{(x-y)^2}dy\geq \frac{u(x)^2}{8\pi\Delta^2} (2\Delta-\frac{4E}{u(x)^2}).$$
We achieve the conclusion of lemma \ref{cuarta} by taking $\Delta=\D\frac{4E}{u(x)^2}$. \quad\eop

\vspace{0.5cm}

Next, let us define $J(t)=u(x(\be_0,t),t)$. Thus, applying  lemma \ref{cuarta} to the expression (\ref{segunda}), we obtain the inequality
\begin{equation}J_{tt}(t)\geq C J(t)^4-J(t).\label{inicio}\end{equation}

Since $Hu_0(\be_0)>0$ and $J_t(t)=Hu(x(\be_0,t),t)$, we obtain that $J_t(t)>0$ and $J(t)>J(0)$ for
$t\in (0,t^*)$ and  $t^*$ small enough. Therefore, multiplying (\ref{inicio}) by $J_t(t)$ we have that
$$\frac{1}{2}(J_t(t)^2)_t\geq \frac{C}{5}(J(t)^5)_t-\frac{1}{2}(J(t)^2)_t,\quad \forall t\in [0,t^*).$$
Integrating this inequality in time from $0$ to $t$ we get
\begin{equation}\label{explota}
J_t(t)\geq \Big(J_t(0)^2+\frac{2C}{5}(J(t)^5-J(0)^5)-(J(t)^2-J(0)^2)\Big)^\frac{1}{2}, \quad \forall t\in [0,t^*).
\end{equation}
Now, since $CJ(0)^4-J(0)\geq 0$, by the statements of the theorem we obtain that $J_{tt}(t)>J_{tt}(0)\geq 0$ for $t\in(0,t^*)$. Therefore, $J_t(t)$ is an increasing function $[0,t^*)$. Thus, the inequality (\ref{explota}) holds for all time $t$ and we have a contradiction.

\begin{rem}
It is easy to check that there exists a large class of functions satisfying the requirement of the theorem (\ref{blowup}). For example, we can consider the function
\begin{eqnarray*}
u_0(x)&=&\frac{-ax}{1+(bx)^2},\\
Hu_0(x)&=&\frac{a}{1+(bx)^2},
\end{eqnarray*}
where $a,\,b>0$. Choosing $a$ and $b$ in a suitable way we can have the norm $||u_0||_{L^2(\ff{R})}$ as small as we want and the norm $||u_0||_{L^\infty(\ff{R})}$ as large as we want.
\end{rem}
\begin{rem}
We note that the requirements (\ref{uno}) and (\ref{dos}) in theorem \ref{blowup} can be replaced by
\begin{eqnarray*}Hu_0(\be_0)&\geq&0,\\u_0(\be_0)&>& \left(32\pi ||u_0||_{L^2(\ff{R})}^2\right)^{1/3},\end{eqnarray*}
attaining the same conclusion.
\end{rem}

\section{Blow up for the whole range $0<\alpha<1$.}
In this  section we shall show formation of singularities for the equation (\ref{BLH}), with $0<\alpha<1$.  The aim is to prove the following result:

\begin{thm}\label{teorema3}
There exist  initial data $u_0 \in L^2(\ff{R})\cap C^{1+\delta}(\ff{R})$, with $0<\delta<1$, and a finite time T, depending on $u_0$, such that
$$\lim_{t\rightarrow T}||u(\cdot,t)||_{C^{1+\delta}(\ff{R})}=\infty$$
where $u(x,t)$ is the solution to the equation (\ref{BLH}).
\end{thm}
Proof: Let us assume that there exists a solution of the equation (\ref{BLH}), $u(x,t)\in C([0,T),C^{1+\delta}(\ff{R}))$, for all time $T<\infty$.
Let $J\pq u$ be the convolution
$$J_q^p u(x)=\int_{\ff{R}}w\pq(x-y)u(y)dy$$
where
\begin{equation*}
w\pq(x)= \left\{\begin{array}{ccc}\frac{1}{|x|^q} \sign(x)&\text{if} \, |x|<1\\ \frac{1}{|x|^p} \sign(x)&\text{if}\, |x|>1\end{array}\right.,
\end{equation*}
with $0<q<1$ and $p>2$.
In order to prove  theorem \ref{teorema3} we shall need the following two lemmas.
\begin{lemma}\label{lambdah}
Let $f$ in $C^{1+\delta}(\ff{R})\cap L^2(\ff{R})$ and $0<\alpha<1$. Then
$$\Lambda^\alpha Hf(x)=k_\al \int_{\ff{R}}\frac{f(x)-f(y)}{|x-y|^{1+\al}}\sign(x-y)dy$$
where $$k_\alpha=-\frac{\Gamma(1+\al)\sin((1+\al)\pi/2)}{\pi}.$$
\end{lemma}

Proof: Let $f$ be a function on the Schwartz class. The inverse Fourier transform formula yields
$$\lae Hf(x)=\frac{1}{2\pi}\int_{\ff{R}}-i\sign(k)|k|^\al \hat{f}(k)\exp(ikx)dk.$$
We will understand the above identity as the following limit
\begin{align*}
\lae Hf(x)
&=\lim_{\ep\to0^+}\frac{1}{2\pi}\int_{\ff{R}}-i\sign(k)|k|^\al \exp(-\ep|k|)\exp(ikx)\Big(\int_{\ff{R}}f(y)\exp(-iky)dy\Big)dk\\
&=\lim_{\ep\to0^+}\frac{1}{\pi}\int_{\ff{R}}f(y)\Big(\int_0^\infty k^\alpha \exp(-\ep k)\sin(k(x-y))dk\Big) dy.
\end{align*}
Next, we can compute that
\begin{align*}
\lae Hf(x)&=\lim_{\ep\to0^+}\frac{\Gamma(1+\al)}{\pi}\int_{\ff{R}}\frac{f(y)}{(\ep^2+(x-y)^2)^{(1+\al)/2}}\sin\Big((1+\al)\arctan\Big(
\frac{x-y}{\ep}\Big)\Big)dy\\
&=-\lim_{\ep\to0^+}\frac{\Gamma(1+\al)}{\pi}\int_{\ff{R}}\frac{f(x)-f(y)}{(\ep^2+(x-y)^2)^{(1+\al)/2}}\sin\Big((1+\al)\arctan\Big(\frac{x-y}{\ep}
\Big)\Big)dy\\
&=-\frac{\Gamma(1+\al)\sin((1+\al)\pi/2)}{\pi}\int_{\ff{R}}\frac{f(x)-f(y)}{|x-y|^{1+\al}}\sign(x-y)dy.
\end{align*}
We achieve the conclusion of lemma \ref{lambdah} by the classical density argument. \quad\eop

\begin{lemma}\label{estimacion}
Let $I\pq(x)$ be the integral
$$I\pq(x)=\int_{\ff{R}}\frac{w\pq(x)-w\pq(y)}{|x-y|^{1+\alpha}}\sign(x-y)dy$$
where $0<q<1$ and $p>2$. Then
$$|I\pq(x)|\leq\left\{ \begin{array}{ccc} \frac{K^1}{|x|^{q+\al}} & \text{if $0<|x|<\frac{1}{2}$}\\ \frac{K^2}{|x|^{2+\alpha}} & \text{ if $2<|x|<\infty$}\\
K^3 & \text{if $\frac{1}{2}\leq |x|\leq 2$}\end{array}\right. $$
where $K^1$, $K^2$ and $K^3$ are universal constants depending on $q$ and $p$.
\end{lemma}

Proof: Since the function $I\pq(x)$ is even, we can assume that $x>0$. The constant values of $K^1$ and $K^2$ can be different along the estimates below.

First, let us consider the case $0<x<1/2$. We split as follows
$$I_q^p(x)=\int_{|y|<1}dy+\int_{|y|>1}dy=I_1(x)+I_2(x).$$
It yields
\begin{align*} I_1(x)&=\int_{|y|<1}\frac{\frac{1}{x^q}-\sign(y)\frac{1}{|y|^q}}{|x-y|^{1+\al}}\sign(x-y)dy\\
&=\int_0^1\Big(\frac{\frac{1}{x^q}-\frac{1}{y^q}}{|x-y|^{1+\al}}\sign(x-y)+\frac{\frac{1}{x^q}+\frac{1}{y^q}}{|x+y|^{1+\al}}\Big)dy,
\end{align*}
and a change of variables allow us to split further
\begin{align*}
I_1(x)&=\frac{1}{x^{q+\al}}\int_0^\frac{1}{x}\Big(\frac{1-\frac{1}{\eta^q}}{|1-\eta|^{1+\alpha}}\sign(1-\eta)+\frac{1+\frac{1}{\eta^q}}{|1+\eta|^{1+\alpha}}\Big)d\eta\\
&=\frac{1}{x^{q+\al}}\Big(\int_{0}^1+\int_1^\frac{1}{x}\Big)=\frac{1}{x^{q+\al}}(F_1(x)+F_2(x)).
\end{align*}
For $F_1(x)$ we find the bound
$$|F_1(x)|\leq \int_0^1\Big|\frac{1-\frac{1}{\eta^q}}{|1-\eta|^{1+\alpha}}\Big|d\eta+\int_0^1\Big|
\frac{1+\frac{1}{\eta^q}}{|1+\eta|^{1+\alpha}}\Big|d\eta\leq K_1.$$
On the other hand
$$F_2(x)=\int_1^\frac{1}{x}\Big(\frac{\frac{1}{\eta^q}-1}{|1-\eta|^{1+\al}}+\frac{\frac{1}{\eta^q}+1}{|1+\eta|^{1+\al}}\Big)d\eta=\int_1^\frac{3}{2}+\int_{\frac{3}{2}}^\frac{1}{x}=j_1(x)+j_2(x).$$
For $j_1(x)$ it is easy to obtain
$$|j_1(x)|\leq \int_1^\frac{3}{2}\Big|\frac{\frac{1}{\eta^q}-1}{|1-\eta|^{1+\al}}\Big|d\eta+\int_1^\frac{3}{2}\Big|\frac{\frac{1}{\eta^q}+1}{|1+\eta|^{1+\al}}\Big|d\eta\leq K_1.$$
For $j_2(x)$ we decompose as follows
$$j_2(x)=\int_{\frac{3}{2}}^\frac{1}{x}\frac{1}{\eta^q}\Big(\frac{1}{|1-\eta|^{1+\alpha}}+\frac{1}{|1+\eta|^{1+\alpha}}\Big)d\eta+\int_\frac{3}{2}^\frac{1}{x}\Big(\frac{1}{|1-\eta|^{1+\alpha}}-\frac{1}{|1+\eta|^{1+\alpha}}\Big)d\eta.$$
Thus, since $0<q<1$ and
$$\Big|\frac{1}{|1-\eta|^{1+\alpha}}+\frac{1}{|1+\eta|^{1+\alpha}}\Big|\leq\frac{ K^1}{|\eta|^{1+\al}},\quad\mbox{for}\quad\eta\in[3/2,\infty) $$
we have that
$$|j_2(x)|\leq K^1\int_{\frac{3}{2}}^\infty\frac{1}{\eta^{q+1}}d\eta+\int_\frac{3}{2}^\infty\Big|\frac{1}{|1-\eta|^{1+\alpha}}-\frac{1}{|1+\eta|^{1+\alpha}}\Big|d\eta\leq K^1.$$
Let us continue with $I_2$ which can be written in the form \begin{align*}
I_2(x)&=\int_{|y|>1}\frac{\frac{1}{x^q}-\sign(y)\frac{1}{|y|^p}}{|x-y|^{1+\al}}\sign(x-y)dy=\int_1^\infty\Big(-\frac{\frac{1}{x^q}-\frac{1}{|y|^p}}{|x-y|^{1+\al}}+\frac{\frac{1}{x^q}+\frac{1}{|y|^p}}{|x+y|^{1+\al}}\Big)dy\\
&=\frac{1}{x^{q+\al}}\int_{\frac{1}{x}}^\infty\Big(\frac{\frac{x^{q-p}}{\eta^p}-1}{|1-\eta|^{1+\al}}+\frac{1+\frac{x^{q-p}}{\eta^p}}{|1-\eta|^{1+\al}}\Big)d\eta.
\end{align*}
The following decomposition
\begin{align*}
I_2(x)&=\frac{1}{x^{q+\al}}\int_{\frac{1}{x}}^\infty\Big(\frac{1}{|1+\eta|^{1+\alpha}}-\frac{1}{|1-\eta|^{1+\alpha}}\Big)d\eta\\
&\quad+\frac{1}{x^{p+\al}}\int_{\frac{1}{x}}^\infty\frac{1}{\eta^p}\Big(\frac{1}{|1-\eta|^{1+\al}}+\frac{1}{|1+\eta|^{1+\alpha}}\Big)d\eta\end{align*}
yields
\begin{align*}
|I_2(x)|&\leq \frac{K^1}{x^{q+\al}}+\frac{1}{x^{p+\al}}\int_{\frac{1}{x}}^\infty \frac{1}{\eta^p}\Big|\frac{1}{|1-\eta|^{1+\al}}+\frac{1}{|1+\eta|^{1+\al}}\Big|d\eta\\
&\leq\frac{K^1}{x^{q+\al}}+\frac{K^1}{x^{p+\al}}\int_{\frac{1}{x}}^\infty \frac{1}{\eta^{p+1}}d\eta\leq K^1\Big(\frac{1}{|x|^{\alpha+q}}+\frac{1}{|x|^{\alpha}}\Big)\leq \frac{K_1}{x^{q+\alpha}}.
\end{align*}

Next, we consider the case $2<x<\infty$ taking
$$I\pq(x)=\int_{\ff{R}}\frac{\frac{1}{x^p}-w(y)}{|x-y|^{1+\al}}\sign(x-y)dy=\int_{|y|<1}dy+\int_{|y|>1}dy=J_1(x)+J_2(x).$$
For $J_2(x)$ we have that
\begin{align*}
J_2(x)&=\int_{|y|>1}\frac{\frac{1}{x^p}-\sign(y)\frac{1}{|y|^p}}{|x-y|^{1+\al}}\sign(x-y)dy\\
&=\int_{1}^{\infty}\Big(\frac{\frac{1}{x^p}-\frac{1}{|y|^p}}{|x-y|^{1+\al}}\sign(x-y)+\frac{\frac{1}{x^p}+\frac{1}{|y|^p}}{|x+y|^{1+\al}}\Big)dy
\end{align*}
and a change of variables provides
\begin{align*}
J_2(x)&=\frac{1}{x^{p+\al}}\int_{\frac{1}{x}}^\infty\Big(\frac{1-\frac{1}{\eta^p}}{|1-\eta|^{1+\alpha}}\sign(1-\eta)+\frac{1+\frac{1}{\eta^p}}{|1+\eta|^{1+\alpha}}\Big)d\eta\\
&=\frac{1}{x^{\al+p}}\Big(\int_{{\frac{1}{x}}}^1+\int_1^\infty\Big)=\frac{1}{x^{\al+p}}(H_1(x)+H_2(x)).
\end{align*}
For $H_2(x)$ one could bound as follow
$$|H_2(x)|\leq \int_{1}^\infty \Big|\frac{1-\frac{1}{\eta^p}}{|1-\eta|^{1+\alpha}}\Big|d\eta+\int_{1}^\infty\Big|\frac{1+\frac{1}{\eta^p}}{|1+\eta|^{1+\alpha}}\Big|d\eta\leq K^2.$$
On the other hand, in $H_1(x)$ we split further
$$H_1(x)=\int_{\frac{1}{x}}^1\Big(\frac{1-\frac{1}{\eta^p}}{|1-\eta|^{1+\alpha}}+\frac{1+\frac{1}{\eta^p}}{|1+\eta|^{1+\alpha}}\Big)d\eta=\int_\frac{1}{x}^\frac{2}{3}d\eta+\int_{\frac{2}{3}}^1d\eta=h_1(x)+h_2(x).$$
The term $h_2(x)$ is bounded by
$$|h_2(x)|\leq \int_{\frac{2}{3}}^1\Big|\frac{1-\frac{1}{\eta^p}}{|1-\eta|^{1+\alpha}}\Big|d\eta+\int_{\frac{2}{3}}^1\Big|\frac{1+\frac{1}{\eta^p}}{|1+\eta|^{1+\alpha}}\Big|d\eta\leq K^2.$$
We reorganize $h_1(x)$ so that
$$h_1(x)=\int_{\frac{1}{x}}^\frac{2}{3}\Big(\frac{1}{|1-\eta|^{1+\al}}+\frac{1}{|1+\eta|^{1+\al}}\Big)d\eta+\int_{\frac{1}{x}}^\frac{2}{3}\frac{1}{\eta^p}\Big(\frac{1}{|1-\eta|^{1+\al}}-\frac{1}{|1+\eta|^{1+\al}}\Big)d\eta.$$
Since $p>2$ and
$$\left|\frac{1}{|1-\eta|^{1+\al}}-\frac{1}{|1+\eta|^{1+\al}}\right|\leq K^2 \eta\quad \mbox{for}\quad \eta\in[0,2/3],$$
we obtain that
$$|h_1(x)|\leq\int_0^\frac{2}{3}\Big|\frac{1}{|1-\eta|^{1+\al}}+\frac{1}{|1+\eta|^{1+\al}}\Big|d\eta+K^2\int_{\frac{1}{x}}^\frac{2}{3}\frac{1}{\eta^{p-1}}d\eta\leq K^2(1+x^{p-2}).$$
Therefore
$$|J_2(x)|\leq K^2\Big(\frac{1}{x^{p+\alpha}}+\frac{1}{x^{2+\alpha}}\Big)\leq \frac{K^2}{x^{2+\al}} .$$
Next, we deal with $J_1$ given by
\begin{align*}
J_1(x)&=\int_{|y|<1}\frac{\frac{1}{x^p}-\sign(y)\frac{1}{|y|^q}}{|x-y|^{1+\alpha}}dy=
\int_0^1\left(\frac{\frac{1}{x^p}-\frac{1}{|y|^q}}{|x-y|^{1+\alpha}}+\frac{\frac{1}{x^p}+\frac{1}{|y|^q}}{|x+y|^{1+\alpha}}\right)dy\\
&=\frac{1}{x^{p+\alpha}}\int_0^\frac{1}{x}\left(\frac{1-\frac{x^{p-q}}{\eta^q}}{|1-\eta|^{1+\alpha}}+\frac{1+\frac{x^{p-q}}{\eta^q}}{|1+\eta|^{1+\alpha}}\right)d\eta.
\end{align*}
Hence
\begin{align*}
J_1(x)&=\frac{1}{x^{p+\alpha}}\!\int_0^\frac{1}{x}\!\Big(\frac{1}{|1-\eta|^{1+\alpha}}+\frac{1}{|1+\eta|^{1+\alpha}}\Big)d\eta+\frac{1}{x^{q+\alpha}}\!\int_0^\frac{1}{x}\!\frac{1}{\eta^q}\Big(\frac{1}{|1+\eta|^{1+\alpha}}-\frac{1}{|1-\eta|^{1+\alpha}}\Big)d\eta.
\end{align*}
Since  $p>2$ and
$$\Big|\frac{1}{|1+\eta|^{1+\alpha}}-\frac{1}{|1-\eta|^{1+\alpha}}\Big|\leq K^2\eta,\quad\mbox{for}\quad \eta\in[0,1/2],$$
we obtain that
$$|J_1(x)|\leq \frac{K^2}{x^{p+\alpha}}+\frac{K^2}{x^{q+\al}}\int_0^\frac{1}{x}\eta^{1-q}d\eta\leq K^2\Big(\frac{1}{x^{p+\alpha}}+\frac{1}{x^{2+\alpha}}\Big)\leq \frac{K^2}{x^{2+\alpha}}.$$
The bound for $1/2<x<2$ is obvious, which allow us to conclude the proof. \quad\eop

\vspace{0.5cm}

In order to prove theorem \ref{teorema3}, we shall study the evolution of $J(t)=J\pq u(x_b(t),t)$, where $x_b$ is the trajectory $x_b(t)=x(0,t)$. Hence
$$\frac{dJ(t)}{dt}=-\frac{1}{2} J\pq ((u^2)_x)(x_b(t),t)+ J\pq\lae H u (x_b(t),t)+u(x_b(t),t)(\pa_x J\pq u)(x_b(t),t).$$
We can write
$$
J\pq ((u^2)_x) = \int_\ff{R}(u(x)^2-u(y)^2)W\pq(x-y) dy
$$
and
$$\partial_x (J\pq u)(x)  = \int_{\ff{R}}(u(x)-u(y))W\pq(x-y)dy$$
where
$$W\pq=\left\{\begin{array}{ccc} \frac{q}{|x|^{q+1}} & \text{if $|x|<1$} \\ \frac{p}{|x|^{p+1}} & \text{if $|x|>1$} \end{array}\right..$$
Then, it is easy to check that
$$-\frac{1}{2} J\pq ((u^2)_x)(x)+u(x)(\pa_x J\pq u)(x)=\frac{1}{2}\int_\ff{R}(u(x)-u(y))^2 W\pq(x-y) dy,$$
and therefore
\begin{equation}\label{pepe}
\frac{dJ(t)}{dt}=\frac{1}{2}\int_\ff{R}(u(x_b(t))-u(y))^2 W\pq(x_b(t)-y) dy+J\pq\lae H u (x_b(t),t).\end{equation}
Using lemma \eqref{lambdah}, the linear term becomes
$$J\pq\lae H u (x)=k_\alpha\int_{\ff{R}}w\pq(x-y)\int_{\ff{R}}\frac{u(y)-u(s)}{|y-s|^{1+\al}}\sign(y-s)dsdy,$$
and a wise use of the principal value provides
\begin{align*}
J\pq\lae H u (x)&=k_\alpha\int_{\ff{R}}w\pq(x-y)P.V.\int_{\ff{R}}\frac{-u(s)}{|y-s|^{1+\al}}\sign(y-s)dsdy\\
&=k_\alpha\int_{\ff{R}}w\pq(x-y)P.V.\int_{\ff{R}}\frac{u(x)-u(s)}{|y-s|^{1+\al}}\sign(y-s)dsdy\\
&=k_\alpha\int_{\ff{R}}(u(x)-u(s))P.V.\int_{\ff{R}}\frac{w\pq(x-y)}{|y-s|^{1+\al}}\sign(y-s)dyds\\
&=k_\alpha\int_{\ff{R}}(u(x)-u(s))\int_{\ff{R}}\frac{w\pq(x-s)-w\pq(x-y)}{|y-s|^{1+\al}}\sign(s-y)dyds
\end{align*}
to find finally
$$
J\pq\lae H u (x)=k_\alpha \int_{\ff{R}}(u(x)-u(s))I\pq(x-s)ds.
$$
Therefore
\begin{align*}
|J\pq\lae H u (x)|&\leq |k_\alpha|\int_{\ff{R}}|u(x)-u(y)||I\pq(x-y)|dy\\
&\leq |k_\alpha|\left(\int_{\ff{R}}(u(x)-u(y))^2 W\pq(x-y)dy\right)^\frac{1}{2}\left(\int_{\ff{R}} \frac{I\pq(x)^2}{W\pq(x)}dx\right)^{\frac{1}{2}}.
\end{align*}
Since
$$\frac{I\pq(x)^2}{W\pq(x)}\leq\left\{\begin{array}{ccc} \frac{C}{|x|^{2\alpha+q-1}} & \text{ when $|x|\rightarrow 0$}\\ \frac{C}{|x|^{3+2\alpha-p}}  & \text{ when $|x|\rightarrow \infty$}\end{array}\right. ,$$
by taking $2<p<2+2\al$ and  $q<2(1-\alpha)$, we obtain that
\begin{align*}
|J\pq\lae H u (x)|&\leq C(p,q)\left(\int_{\ff{R}}(u(x)-u(y))^2 W\pq(y)dy\right)^\frac{1}{2}\\
&\leq\frac{1}{4}\int_{\ff{R}}(u(x)-u(y))^2 W\pq(x-y)dy+ C.
\end{align*}
This inequality in the equation (\ref{pepe}) yields
$$\frac{dJ(t)}{dt}\geq\frac{1}{4}\int_\ff{R}(u(x_b(t))-u(y))^2 W\pq(x_b(t)-y) dy-C(p,q)$$
On the other hand
$$J(t)=\int_{\ff{R}}u(y)w\pq(x_b(t)-y)dy=\int_{\ff{R}}(u(y)-u(x_b(t)))w\pq(x_b(t)-y)dy$$
$$\leq\left(\int_{\ff{R}}(u(x_b(t))-u(y))^2 W\pq (x_b(t)-y)dy\right)^\frac{1}{2}\left(\int_{\ff{R}}\frac{w\pq(x)^2}{W\pq(x)}dx\right)^\frac{1}{2}.$$
The following bound
$$\frac{w\pq(x)^2}{W\pq(x)}\leq\left\{\begin{array}{ccc} \frac{C}{|x|^{q-1}} & \text{ when $|x|\rightarrow 0$}\\ \frac{C}{|x|^{p-1}}  & \text{ when $|x|\rightarrow \infty$}\end{array}\right.$$
allows us to obtain, for $2<p<2+2\al$ and $0<q<1$,
$$\int_{\ff{R}}(u(x_b(t))-u(y))^2 W\pq(x_b(t)-y)dy\geq c(q,p)J(t)^2.$$
Therefore we obtain a quadratic evolution equation
$$\frac{dJ(t)}{dt}\geq c(q,p) J(t)^2-C(q,p)$$
and by taking $c(q,p) J(0)^2-C(q,p)>0$, we find a contradiction for the mere fact that
$$J(t)\leq C(q,p)||u||_{L^\infty}.$$

\section{Appendix}

In this section we generalize the pointwise inequality \eqref{pwi} evolving the nonlocal operator $2f\Lambda^\al f-\Lambda^\al(f^2).$ Some simple applications to Gagliardo-Nirenberg-Sobolev inequalities are also shown.

\begin{lemma}
Consider a function $f:\R^n\to\R$ in the Schwartz class, $0<\al<2$ and $0<p<\infty$. Then
\begin{equation}\label{pwig}
|f(x)|^{2+\frac{\al p}{n}}\leq C(\al,p,n)||f||_{L^p(\ff{R}^n)}^{\frac{\al p}{n}}(2f(x)\Lambda^\al f(x)-\Lambda^\al(f^2)(x))
\end{equation}
for any $x\in\R^n$.
\end{lemma}

Proof: The formula for the operator $\Lambda^\al$ in $\R^n$
$$\Lambda^\al f(x)=k_{\alpha,n}\int_{\R^n}\frac{f(x)-f(y)}{|x-y|^{n+\al}}dy$$
and $0<\al<2$, allows us to find
$$
2f(x)\Lambda^\al f(x)-\Lambda^\al(f^2)(x)=k_{\al,n}\int_{\ff{R}^n}\frac{(f(y)-f(x))^2}{|x-y|^{n+\al}}dy.
$$
We consider $f(x)>0$, being the case $f(x)<0$ analogous. Let $\Omega$, $\Omega^1$ and $\Omega^2$ be the sets
\begin{align*}
\Omega\,\,\,&=\{y\in \ff{R}\,:|x-y|<\Delta\},\\
\Omega^1&=\{y\in \Omega\, : f(x)-f(y)\geq f(x)/2\},\\
\Omega^2&=\{y\in \Omega\, : f(x)-f(y)< f(x)/2\}=\{y\in \Omega \, : f(y)> f(x)/2\}.
\end{align*}
Then
\begin{equation*}
2f(x)\Lambda^\al f(x)-\Lambda^\al(f^2)(x)\geq k_{\al,n}\int_{\Omega} \frac{(f(y)-f(x))^2}{|x-y|^{n+\al}}dy \geq  k_{\al,n}\frac{f(x)^2}{4 \Delta^{n+\al}}|\Omega^1|.
\end{equation*}
On the other hand
\begin{equation*}
||f||_{L^p(\R^n)}^p=\int_{\R^n}|f(y)|^p dy\geq \frac{f(x)^p}{2^p}|\Omega^2|,
\end{equation*}
therefore \begin{equation*}
2f(x)\Lambda^\al f(x)-\Lambda^\al(f^2)(x)\geq k_{\al,n} \frac{f(x)^2}{4\Delta^{n+\alpha}}(|\om|-|\om^2|)\geq k_{\al,n} \frac{f(x)^2}{4\Delta^{n+\alpha}}(c_n\Delta^n-\frac{2^p||f||_{L^p(\ff{R}^n)}^p}{f(x)^p}),\end{equation*}
where $c_n=2\pi^{n/2}/(n\Gamma(n/2))$. By choosing
\begin{equation*}
\Delta^n=\frac{(n+\al)2^p||f||^p_{L^p(\ff{R}^n)}}{\al c_n f(x)^p}
\end{equation*}
we obtain the desired estimate.  \quad\eop

\begin{rem}
Inequality \eqref{pwig} allows us to get easily the following Gagliardo-Nirenberg-Sobolev estimate:
\begin{equation*}
||f||_{L^{2+\frac{\al p}{n}}}^{2+\frac{\al p}{n}}\leq 2C(\al,p,n)||f||_{L^p(\ff{R}^n)}^{\frac{\al p}{n}}||\Lambda^\frac{\alpha}{2} f||_{L^2(\ff{R}^n)}^2,
\end{equation*}
for $0<\al<2$ and $0<p<\infty$.
\end{rem}

\subsection*{{\bf Acknowledgments}}

\smallskip

The authors wish to thank John Hunter and Peter Constantin for helpful discussions, and Edriss Titi for suggesting the problem.

The authors were partially supported by the grant {\sc MTM2008-03754} of the MCINN (Spain) and
the grant StG-203138CDSIF of the ERC. The third author was partially supported by NSF-DMS grant
0901810.

\vspace{1cm}

\begin{quote}
\begin{tabular}{l}
\textbf{Angel Castro} \\
{\small Instituto de Ciencias Matem\'aticas}\\
{\small Consejo Superior de Investigaciones Cient\'ificas}\\
{\small Serrano 123, 28006 Madrid, Spain}\\
{\small Email: angel\underline{  }castro@icmat.es}
\end{tabular}
\end{quote}
\begin{quote}
\begin{tabular}{ll}
\textbf{Diego C\'ordoba} &  \textbf{Francisco Gancedo}\\
{\small Instituto de Ciencias Matem\'aticas} & {\small Department of Mathematics}\\
{\small Consejo Superior de Investigaciones Cient\'ificas} & {\small University of Chicago}\\
{\small Serrano 123, 28006 Madrid, Spain} & {\small 5734 University Avenue, Chicago, IL 60637}\\
{\small Email: dcg@icmat.es} & {\small Email: fgancedo@math.uchicago.edu}
\end{tabular}
\end{quote}

\end{document}